\documentclass{birkjour}
\usepackage[noadjust]{cite}
\usepackage{xcolor}
\usepackage{longtable}
\usepackage{wrapfig}
\usepackage{url}
\newtheorem{thm}{Theorem}[section]
\newtheorem{cor}[thm]{Corollary}
\newtheorem{lem}[thm]{Lemma}
\newtheorem{prop}[thm]{Proposition}
\theoremstyle{definition}
\newtheorem{defn}[thm]{Definition}
\theoremstyle{remark}

\numberwithin{equation}{section}

\newcommand{\BibTeX}{\textsc{B\kern-0.1emi\kern-0.017emb}\kern-0.15em\TeX}
\newcommand{\ed}{\end{document}}
\newcommand{\tabincell}[2]{\begin{tabular}{@{}#1@{}}#2\end{tabular}}

\begin{document}

%
%
%
%
%
%
%
%
%
\pdfoutput=1
\title[Generalized derivations of Complex $\omega$-Lie
Superalgebras]
 {Generalized Derivations\\ of Complex $\omega$-Lie
Superalgebras}
\author[Jia zhou]{Jia Zhou}
%

\address{%
College of Information Technology, Jilin Agricultural University\\
2888, Xincheng Str.\\
 Nanguan District, Changchun, China}
\email{zhoujia@jlau.edu.cn}

\date{\today}
\begin{abstract}
~Let
$(g,~[-,-],~\omega)$ be a finite-dimensional complex $\omega$-Lie superalgebra. This paper explores the algbaraic structures of generalized derivation superalgebra ${\rm GDer}(g)$, compatatible generalized derivations algebra ${\rm GDer}^{\omega}(g)$, and their subvarieties such as quasiderivation superalgebra ${\rm QDer}(g)$(${\rm QDer}^{\omega}(g)$), centroid ${\rm Cent}(g)$ (${\rm Cent}^{\omega}(g)$) and quasicentroid ${\rm QCent}(g)$
(${\rm QCent}^{\omega}(g)$). We prove that ${\rm GDer}^{\omega}(g) = {\rm
QDer}^{\omega}(g) + {\rm QCent}^{\omega}(g)$. We also study the embedding question of compatible quasiderivations of $\omega$-Lie superalgebras, demonstrating that ${\rm QDer}^{\omega}(g)$ can be embedded as derivations in a larger $\omega$-Lie superalgebra $\breve g$ and furthermore, we obtain a semidirect sum decomposition: ${\rm Der}^{\omega}(\breve{g})=\varphi({\rm
QDer}^{\omega}(g))\oplus {\rm ZDer}(\breve{g})$, when the annihilator of $g$ is zero. In particular, for the 3-dimensional complex $\omega$-Lie superalgebra $H$, we explicitly calculate ${\rm GDer}(H)$, ${\rm GDer}^{\omega}(H)$, ${\rm QDer}(H)$ and ${\rm QDer}^{\omega}(H)$, and derive the Jordan standard forms of generic elements in these varieties.
\end{abstract}

\maketitle
\section{Introduction}

The notion of
$\omega$-Lie algebras, which is related to the study of
isoparametric hypersufaces in Riemannian geometry (See Bobienski-Nurowski
\cite{Bobienski1} and Nurowski \cite{Nurowski2}), was initially introduced by Nurowski \cite{Nurowski3}, who gave the first example of nontrivial 3-dimensional $\omega$-Lie algebra and a classification of 3-dimensional $\omega$-Lie algebras over real numbers. Several conclusions about algebraic structures of $\omega$-Lie algebras have been gradually achieved for the past 10 years. In 2014, Chen-Liu-Zhang \cite{Chen1} obtained a classification of 3-dimensional complex $\omega$-Lie algebras. With the list of 3-dimensional classification, Chen-Zhang \cite{Chen2} completed a classification of 4-dimensional complex Lie algebras. In particular, they obtained a classification of nontrivial finite-dimensional complex simple $\omega$-Lie algebras, see (\cite{Chen2}, Theorem 1.7). Among those studies, derivation theory plays an essential role in understanding structures, representations, extensions and automorphism groups of $\omega$-Lie algebras.  Chen-Zhang \cite{Chen3} show that the set ${\rm Der}^{\omega}(g)$ of all $\omega$-derivations is a Lie subalgebra of ${\rm Der}(g)$ and the set${\rm Aut}^{\omega}(g)$ of
all $\omega$-automorphisms is a subgroup of ${\rm Aut}(g)$. They also compute explicitly ${\rm Der}(g)$ and ${\rm Aut}(g)$ for any 3-dimensional and 4-dimensional nontrivial $\omega$-Lie algebra g, show that all 3-dimensional nontrivial $\omega$-Lie algebras are multiplicative, and provide a 4-dimensional example of $\omega$-Lie algebra that is not multiplicative. For other statements see \cite{Zhang1,Hassan}.

As an important subject, Generalized derivation theory of Lie algebras was firstly introduced by Leger and Luks. In their paper \cite{Leger}, they obtained some properties of the generalized derivation algebras and their subalgebras such as quasiderivation algebras and centroids, they proved that a Lie algebra $g$ satisfying ${\rm QDer}(g)={\rm End}(g)$ is either abelian, two-dimensional solvable or three-dimensional simple. In particular, they characterized the relationship between quasiderivations and cohomology of Lie algebras. Zhang-Zhang \cite{Zhang1} generalized some results in \cite{Leger} to the generalized derivation superalgebras of a Lie superalgebra $L$. They studied the derivation superalgebra ${\rm Der}(L),$ the quasiderivation superalgebras ${\rm QDer}(L)$ and the generalized derivation superalgebras ${\rm GDer}(L),$ which formed a tower ${\rm Der}(L)\subseteq {\rm QDer}(L)\subseteq {\rm GDer}(L)\subseteq {\rm End}(L).$ More precisely, they characterized the Lie superalgebra for which the generalized derivation superalgebras or their Lie subsuperalgebras satisfy some special conditions. Furthermore, they obtained a semidirect sum decomposition of ${\rm Der}(\tilde L)$ of a larger Lie subsuperalgebra $\tilde L$, when the annihilator of $L$ is zero. Article \cite{Zhou2} contains some properties of the generalized derivation algebra ${\rm
GDer}(T)$ of a Hom-Lie triple system $T$, a equation for the quasiderivation algebra and the
quasicentroid as ${\rm GDer}(T)={\rm QDer}(T)+{\rm QC}(T)$, a conclusion about quasiderivation algebra ${\rm QDer}(T)$ and derivations of
a larger Hom-Lie triple system $\tilde T$. General results on centroids of a Hom-Lie triple system $T$
are also developed in this paper. For the generalized derivation algebras of a variety of other non-associative algebras, the readers will be referred to \cite{Zhang2,Zhou1,Zhou3,Zhou4,Liangyun1,Ataguema,Bai,Kaygorodov}.

The notation of $\omega$-Lie superalgebras was introduced by Zhou-Chen-Ma in \cite{Zhou5}. Let $\mathbb K$ be a field of characteristic zero. An $\omega$-Lie superalgebra $g=g_{\bar{0}}\oplus g_{\bar{1}}$ is a $\mathbb Z_2$-graded vector space over $\mathbb K$ with a multiplication $[\cdot,\cdot]:g\times g\rightarrow g$ and a bilinear form $\omega: g\times g\rightarrow \mathbb K$ satisfying

 $(1)$\, $[g_{\alpha},g_{\beta}]\subseteq g_{\alpha+\beta},\quad \forall~\alpha,\beta \in \mathbb Z_2,$

 $(2)$\, $[x,y]=-(-1)^{\mid x\mid \mid y\mid}[y,x],$~~~~~~~~~~~~~(graded skew-symmetric)

 $(3)$\, $(-1)^{\mid x| | y|}[[y,z],x]+(-1)^{| y| | z|}[[z,x],y]+(-1)^{| x| | z|}[[x,y],z]=\\(-1)^{| x| | y|}\omega(y,z)x+(-1)^{| y| | z|}\omega(z,x)y+(-1)^{| x| | z|}\omega(x,y)z,$ ~~~~~~(graded $\omega$-Jacobi identity)

$(4)$\, $\omega(g_{\bar{0}},g_{\bar{1}})=0,$ for all homogeneous elements $x,y,z\in g$.
As a $\mathbb Z_2$-graded vector space,
an $\omega$-Lie superalgebra $g=g_{\bar{0}}\oplus g_{\bar{1}}$ is an $\omega$-Lie algebra
if $g_{\bar{1}}=0$ i.e., $g=g_{\bar{0}}$. An $\omega$-Lie
superalgebra, with the bilinear form
$\omega,$ becomes a Lie superalgebra if $\omega\equiv0.$ Hence we usually call Lie superalgebras
trivial $\omega$-Lie superalgebras. Therefore, $\omega$-Lie superalgebras generalize both $\omega$-Lie algebras and Lie superalgebras.  Some elementary properties and representations of $\omega$-Lie superalgebras were characterized in the literature \cite{Zhou5}, which contained all classifications for 3-dimensional and 4-dimensional $\omega$-Lie superalgebras over the field of complex numbers. In \cite{Zhou6}, Zhou-Chen-Ma introduce the notions of derivation superalgebra ${\rm Der}(g)$, superalgebra ${\rm Der}^{\omega}(g)$, automorphism group ${\rm Aut}(g)$ and subgroup of ${\rm Aut}^{\omega}(g)$ of a $\omega$-Lie algebra $g$. For any 3-dimensional or 4-dimensional  nontrivial non-$\omega$-Lie complex $\omega$-Lie superalgebra $g$, they calculate explicitly ${\rm Der}(g)$ and ${\rm Aut}(g)$, give a conclusion that $g$ is multiplicative, as well as they show that any irreducible respresentation of 4-dimensional $\omega$-Lie superalgebra $P_{2,k}(k\ne 0,-1)$ is 1-dimensional.

The purpose of this article is to develop a theory of generalized derivations for finite-dimensional $\omega$-Lie superalgebras. We proceed as follows.
Section 2 contains some basic definitions which will be used in what follows. In section 3, we denote the set of all compatible linear maps $d$ in ${\rm Der}(g),$ ${\rm QDer}(g),$ ${\rm C}(g),$ ${\rm QC}(g),$ ${\rm ZDer}(g)$ by ${\rm Der}^{\omega}(g),$ ${\rm QDer}^{\omega}(g),$ ${\rm C}^{\omega}(g),$ ${\rm QC}^{\omega}(g)$ and ${\rm ZDer}^{\omega}(g)$, respectively. We give some fundamental properties about these sets, and a tower ${\rm Der}^{\omega}(g)\subseteq {\rm
QDer}^{\omega}(g)\subseteq {\rm GDer}^{\omega}(g)\subseteq {\rm GDer}(g)\subseteq {\rm
End}(g)$ is formed. Especially, we obtain relationships among ${\rm GDer}^{\omega}(g),$ ${\rm QDer}^{\omega}(g),$ ${\rm C}^{\omega}(g),$ and ${\rm QC}^{\omega}(g)$, see Proposition 3.7 and Theorem 3.10.
In Section 4, we prove that the
quasiderivations of an $\omega$-Lie superalgebra $g$ can be embedded as derivations in a larger
$\omega$-Lie superalgebra $\breve{g}$ and obtain a direct sum
decomposition of ${\rm Der}^{\omega}(\breve{g})$ when the annihilator of $g$
is equal to zero. Section 5 is devoted to calculating for the only nontrivial non-$\omega$-Lie 3-dimensional complex $\omega$-Lie superalgebra $H$, we give the detailed arguments about the Jordan standard forms of ${\rm GDer}(H)$, ${\rm GDer}^{\omega}(H)$, ${\rm QDer}(H)$ and ${\rm QDer}^{\omega}(H)$, and list them in the tables.

Throughout this paper we assume that the ground field is $\mathbb K$ of characteristic zero. An element $x\in g$ is called homogeneous if $x$ is in $g_{\bar{0}}$ or $g_{\bar{1}}$. For any homogeneous element $x$ we shall use the standard notation $|x|\in {\mathbb Z}_2=\{\bar{0},\bar{1}\}$ to indicate its degree.

\section{Basic concepts and properties}

\begin{defn}\cite{Zhou6}
Let $g$ be an $\omega$-Lie
superalgebra.  An homogeneous linear map $d\in {\rm End}(g)$ is said to be a derivation of degree $\mid d\mid$ if
\begin{align*}d([x,y]=[d(x),y]+(-1)^{\mid d\mid\mid x\mid}[x,d(y)],~\forall x,y\in g.\end{align*}

We denote the set of all even derivations and odd derivations by ${\rm Der}_{\bar0}(g)$ and ${\rm Der}_{\bar1}(g)$, respectively. ${\rm
Der}(g)={\rm Der}_{\bar{0}}(g)\bigoplus {\rm Der}_{\bar{1}}(g)$ provided with the Lie-super commutator is a subalgebra of
${\rm End}(g)$ and is called the derivation algebra of $g$.
\end{defn}

\begin{defn}
{Let $g$ be an $\omega$-Lie
superalgebra.  A homogeneous linear map $d\in {\rm End}(g)$ is said to be a generalized derivation of degree $\mid d\mid$ if there exist two homogeneous linear maps $d'$, $d''$ with the same degree of $d$, such that
\begin{align*}d''([x,y]=[d(x),y]+(-1)^{\mid d\mid\mid x\mid}[x,d'(y)],~\forall x,y\in g.\end{align*}}
\end{defn}

\begin{defn}
{Let $g$ be an $\omega$-Lie
superalgebra.  An homogeneous linear map $d\in {\rm End}(g)$ is said to be a homogeneous quasi-derivation of degree $\mid d\mid$ if there exist an endomorphism $d'$ with the same degree of $d$, such that
\begin{align*}d'([x,y]=[d(x),y]+(-1)^{\mid d\mid\mid x\mid}[x,d(y)],~\forall x,y\in g.\end{align*}}
\end{defn}

Let ${\rm GDer}(g)$ and ${\rm QDer}(g)$ be
the sets of generalized derivations and of
quasi-derivations, respectively.$${\rm
GDer}(g)={\rm GDer}_{\bar{0}}(g)\bigoplus {\rm GDer}_{\bar{1}}(g),$$ and $${\rm
QDer}(g)={\rm QDer}_{\bar{0}}(g)\bigoplus {\rm QDer}_{\bar{1}}(g).$$
 It is easy to verify that both ${\rm GDer}(g)$ and
${\rm QDer}(g)$ provided with the Lie-super commutator are subalgebras of
${\rm End}(g).$

\begin{defn}
{Let $g$ be an $\omega$-Lie
superalgebra. ${\rm C}(g)$ is called centroid of $g$, if it is the set consisting of all homogeneous linear maps $d$ of $g$ such that
\begin{align*}
[d(x),y]=(-1)^{\mid d\mid\mid x\mid}[x,d(y)]=d([x,y]),~\forall x,y\in g.\end{align*}is called centroid of $g$.}
\end{defn}

\begin{defn}
{Let $g$ be an $\omega$-Lie
superalgebra.  The set ${\rm QC}(g)$ consists of all homogeneous linear maps $d$ of $g$ such that
\begin{align*}
[d(x),y]=(-1)^{\mid d\mid\mid x\mid}[x,d(y)],~\forall x,y\in g.\end{align*}is called quasi-centroid of $g$. }
\end{defn}

\begin{defn}
{Let $g$ be an $\omega$-Lie
superalgebra. ${\rm
Z}(g)={\rm Z}_{\bar{0}}(g)\bigoplus {\rm Z}_{\bar{1}}(g)=\{x\in g: [x,y]=0, \forall y\in g\}$ is called the center of $g$. A homogeneous linear map $d$ is called a center derivation if $d$ satisfies $$[d(x),y]=d([x,y])=0,$$ for all
$x,~y\in g.$ We use ${\rm ZDer}(g)$ to denote the set of all center derivations.}
\end{defn}

It is easy to verify that $${\rm ZDer}(g)\subseteq {\rm
Der}(g)\subseteq {\rm QDer}(g)\subseteq {\rm GDer}(g)\subseteq {\rm
End}(g).$$
$${\rm C}(g)\subseteq {\rm QC}(g)\subseteq {\rm QDer}(g).$$

\section{Generalized derivation algebras and their subalgebras}

This section contains some basic properties
of center derivation algebra, quasi-derivation algebra and the
generalized derivation algebra of an $\omega$-Lie superalgebra.

\begin{defn}
{Let $g$ be an $\omega$-Lie
superalgebra. A homogeneous linear map $d$ of $g$ is compatible if \begin{align*}\omega(d(x),y)+(-1)^{\mid d\mid\mid x\mid}\omega(x,d(y))=0,~\forall x,y\in g.\end{align*}}
\end{defn}

We write ${\rm GDer}^{\omega}(g)$ for the set of all compatible generalized derivations $d$ in ${\rm GDer}(g)$, which is a vector space satisfying ${\rm GDer}^{\omega}(g)\subseteq {\rm GDer}(g)$. Similarly, all compatible linear maps $d$ in ${\rm Der}(g),$ ${\rm QDer}(g),$ ${\rm C}(g),$ ${\rm QC}(g),$ ${\rm ZDer}(g)$ form subsets ${\rm Der}^{\omega}(g),$ ${\rm QDer}^{\omega}(g),$ ${\rm C}^{\omega}(g),$ ${\rm QC}^{\omega}(g),$ ${\rm ZDer}^{\omega}(g)$, respectively.

Clearly, we have the tower
$${\rm Der}^{\omega}(g)\subseteq {\rm
QDer}^{\omega}(g)\subseteq {\rm GDer}^{\omega}(g)\subseteq {\rm GDer}(g)\subseteq {\rm
End}(g).$$

\begin{lem}
{\it $d,d'$ are compatible homogeneous linear maps, they have the same degrees, then $kd$, $d+d'$ are also compatible with degree $\mid d\mid$.}
\end{lem}

\begin{lem}
{\it $d,d'$ are compatible homogeneous linear maps, then $[d,d']$ is also compatible with degree $\mid d\mid+\mid d'\mid$.}
\end{lem}

\begin{proof}
Indeed, for any homogeneous element $x\in g$, $y\in g,$ we have
\begin{eqnarray*}
  \omega([d,d'](x),y)&=&\omega(d\cdot d'(x)-(-1)^{\mid d\mid\mid d'\mid}d'\cdot d(x),y)\\
   &=&\omega(d\cdot d'(x),y)-(-1)^{\mid d\mid\mid d'\mid}\omega(d'\cdot d(x),y)\\
   &=&-(-1)^{\mid d\mid(\mid d'\mid+\mid x\mid)}\omega(d'(x),d(y))+(-1)^{\mid d'\mid \mid x\mid
   }\omega(d(x),d'(y)).
\end{eqnarray*}
On the other hand,
\begin{eqnarray*}
  \omega(x,[d,d'](y))&=&\omega(x,d\cdot d'(y)-(-1)^{\mid d\mid\mid d'\mid}d'\cdot d(y))\\
   &=&\omega(x,d\cdot d'(y))-(-1)^{\mid d\mid \mid d'\mid}\omega(x,d'\cdot d(y))\\
   &=&-(-1)^{\mid d\mid\mid x\mid}\omega(d(x),d'(y))+(-1)^{\mid d'\mid(\mid x\mid +\mid d\mid)
   }\omega(d'(x),d(y)).
\end{eqnarray*}
Thus $$\omega([d,d'](x),y)+ (-1)^{(\mid d\mid+\mid d'\mid)\mid x\mid}\omega(x,[d,d'](y))=0.$$
\end{proof}

\begin{prop}
{\it Let $g$ be an $\omega$-Lie superalgebra. Then
the following statements hold:

$(1)$ ${\rm GDer}^{\omega}(g)$, ${\rm QGer}^{\omega}(g)$, ${\rm Der}^{\omega}(g)$ and ${\rm C}^{\omega}(g)$ are
subalgebras of ${\rm GDer}(g)$, ${\rm QGer}(g)$, ${\rm Der}(g)$ and ${\rm C}(g)$, respectively.

$(2)$ ${\rm ZDer}^{\omega}(g)$ is an ideal of ${\rm Der}^{\omega}(g)$.}
\end{prop}

\begin{proof}
$(1)$ For any homogeneous linear map $d,d'\in {\rm GDer}^{\omega}(g)\subseteq {\rm GDer}(g),$ from the conclusions in {\cite{Zhang2}, Proposition 2.1}, $[d,d']\in {\rm GDer}(g)$ with degree $\mid d\mid+\mid d'\mid$. By Lemma 3.3, $[d,d']$ is a compatible generalized derivation of $g$, which means ${\rm GDer}^{\omega}(g)$ is a subalgebra ${\rm GDer}(g)$.
Other statements hold by applying the same arguments.

$(2)$ Assume that $d\in {\rm ZDer}^{\omega}(g)$ with degree $\mid d\mid$ and $d'\in {\rm Der}^{\omega}(g)$ with degree $\mid d'\mid$, for any homogeneous $x\in g$, $\forall y\in g.$ Note that
$$\begin{array}{ll}[d,d']([x,y])&=dd'([x,y])-(-1)^{\mid d\mid\mid d'\mid}d'd([x,y])\\
&=d([d'(x),y]+(-1)^{\mid d'\mid\mid x\mid}[x,d'(x)])=0.\end{array}$$
and$$\begin{array}{ll}[[d,d'](x),y]&=[(dd'-(-1)^{\mid d\mid\mid d'\mid}d'd)(x),y]\\&=-(-1)^{\mid d\mid
\mid d'\mid}(d'([d(x),y]-(-1)^{\mid d'\mid(d+\mid x\mid)}[d(x),d'(y)])\\&=-(-1)^{\mid d'\mid(\mid d\mid+\mid x\mid)}[d(x),d'(y)]=0.\end{array}$$
Then $[d,d']\in {\rm ZDer}(g)$ of
degree $\mid d\mid+\mid d'\mid$. By Lemma 3.3, $[d,d']$ is a compatible map of $g$, so ${\rm
ZDer}^{\omega}(g)$ is an ideal of ${\rm Der}^{\omega}(g)$.
\end{proof}

\begin{prop}
${\rm ZDer}(g)={\rm ZGer}^{\omega}(g)$

\end{prop}
\begin{proof}
Assume that $x,y,z\in g$ are linearly independent. If $x\in Z(g),$ then from the grading-$\omega-$identity with $x,y,z$, we have that $\omega(x,y)=0,$ for all $y\in g.$

For all homogeneous $d\in {\rm ZDer}(g)$, note that $d(g)\subseteq {\rm Z}(g),$  it follows that for all homogeneous element $y\in g$, $z\in g$, $\omega(d(y),z)+(-1)^{\mid d\mid\mid y\mid}\omega(y,d(z))=0,$  which means $d$ is compatible i.e. ${\rm ZDer}(g)\subseteq{\rm ZGer}^{\omega}(g)$.
\end{proof}

\begin{thm}Let {\it $g$ be an $\omega$-Lie
superalgebra, then

    $(1)$ $[{\rm Der}^{\omega}(g),{\rm C}^{\omega}(g)]\subseteq {\rm C}^{\omega}(g).$

    $(2)$ $[{\rm QDer}^{\omega}(g),{\rm QC}^{\omega}(g)]\subseteq {\rm QC}^{\omega}(g).$

    $(3)$ $[{\rm QC}^{\omega}(g),{\rm QC}^{\omega}(g)]\subseteq {\rm QDer}^{\omega}(g).$

    $(4)$ ${\rm C}^{\omega}(g)\subseteq {\rm QDer}^{\omega}(g).$}
\end{thm}
\begin{proof}
Firstly, we come to show statement $(1)$, and the proof of $(2)$ are similar. For any $d\in {\rm Der}^{\omega}(g)$ with degree $\mid d\mid$ and $d'\in {\rm C}^{\omega}(g)$ with degree $\mid d'\mid$, similar arguments with Lemma $2.3$ in (\cite{Zhang2}, Chapter 2) show that $[d,d']\in {\rm C}(g)$. Applying Lemma $3.3$ we get that $[d,d']$ is a compatible map of $g$, so $[d,d']\in {\rm C}^{\omega}(g)$ i.e., $[{\rm Der}^{\omega}(g),{\rm C}^{\omega}(g)]\subseteq {\rm C}^{\omega}(g).$

 To prove statement $(3)$, we assume that $d\in{\rm
QC}^{\omega}(g),d'\in{\rm QC}^{\omega}(g)$, homogenous element $x\in
g$ and $y\in g.$
From graded skew-symmetry $[d,d']=-(-1)^{\mid d\mid\mid d'\mid}[d',d]$ and direct calculations, we have
$$[[d,d'](x),y]+(-1)^{(\mid d\mid+\mid d'\mid)x}[x,[d,d'](y)]=0.$$
Let $d''=0,$ we combining with Lemma 3.3, we obtain that $[d,d']\in {\rm
QDer}^{\omega}(g),$ of degree $\mid d\mid+\mid d'\mid$ as desired.

The statement $(4)$ holds that if $d\in {\rm C}^{\omega}(g)\subseteq {\rm C}(g)$, then $d\in{\rm QDer}(g)$ by Lemma $2.3$ in (\cite{Zhang2}, Chapter 2). It follow from Lemma $3.3$ that $d$ is compatible, so ${\rm C}^{\omega}(g)\subseteq {\rm QDer}^{\omega}(g).$
\end{proof}

\begin{prop}
{\it Let $g$ be
an ${\omega}$-Lie superalgebra. If ${\rm QDer}^{\omega}(g)={\rm QDer}(g)$ or ${\rm QC}^{\omega}(g)={\rm QC}(g)$, then
$${\rm GDer}^{\omega}(g)={\rm
QDer}^{\omega}(g)+{\rm QC}^{\omega}(g).$$}
\end{prop}

\begin{proof}
To prove the first side ${\rm GDer}^{\omega}(g)\subseteq{\rm QDer}^{\omega}(g)+{\rm QC}^{\omega}(g),$ we take
$d\in {\rm GDer}^{\omega}(g)$ with degree $\mid d\mid$,
then for all homogeneous elements $x,y\in g$, there exist
 homogenous linear maps $d',d''$ with the same degree with $d$ such that
$$[d(x),y]+(-1)^{\mid d\mid \mid x\mid}[x,d'(y)]=d''([x,y]).$$
By the graded skew-symmetry of $\omega-$ Lie superalgebras, one gets that
$$(-1)^{(\mid d\mid+\mid x\mid)\mid y\mid}[y,d(x)]+(-1)^{\mid x\mid+\mid y\mid}[d'(y),x]=(-1)^{\mid x\mid+\mid y\mid}d''([y,x]),$$
which means
$$[d'(y),x]+(-1)^{\mid d\mid\mid y\mid}[y,d(x)]=d''([y,x]).$$
Hence $d'\in {\rm GDer}(g).$ From equations above, we have
$$[\frac{d+d'}{2}(x),y]+(-1)^{\mid d\mid\mid x\mid}[x,\frac{d+d'}{2}(y)]=d''([x,y]),$$
and
$$[\frac{d-d'}{2}(x),y]-(-1)^{\mid d\mid \mid x\mid}[x,\frac{d-d'}{2}(y)]=0,$$
which imply that $\frac{d+d'}{2}\in {\rm
QDer}(g)$ and $\frac{d-d'}{2}\in {\rm
QC}(g)$, of degree $\mid d\mid$.

Hence$$d=\frac{d+d'}{2}+\frac{d-d'}{2}\in
{\rm QDer}(g)+{\rm QC}(g),$$
i.e., $${\rm GDer}^{\omega}(g)\subseteq{\rm QDer}(g)+{\rm QC}(g).$$
 As ${\rm QDer}^{\omega}(g)={\rm QDer}(g)$ or ${\rm QC}^{\omega}(g)={\rm QC}(g)$, it follows that either $d=\frac{d+d'}{2}$ or $\frac{d-d'}{2}$ is
 compatible. By Lemma 3.2, we get that both $d=\frac{d+d'}{2}$ and $\frac{d-d'}{2}$ are compatible. Hence, $${\rm GDer}^{\omega}(g)\subseteq{\rm QDer}^{\omega}(g)+{\rm QC}^{\omega}(g).$$

To prove the other side ${\rm QDer}^{\omega}(g)+{\rm QC}^{\omega}(g)\subseteq {\rm GDer}^{\omega}(g),$  we recall Proposition $2.4$ in (\cite{Zhang2}, Chapter 2) and get ${\rm QDer}(g)+{\rm QC}(g)\subseteq{\rm GDer}(g),$ which asserts that ${\rm QDer}^{\omega}(g)+{\rm QC}^{\omega}(g)\subseteq{\rm GDer}(g)$. By Lemma 3.2, we see that
$${\rm QDer}^{\omega}(g)+{\rm QC}^{\omega}(g)\subseteq{\rm GDer}^{\omega}(g).$$
\end{proof}

\begin{cor}
\it
Let $g$ be
an ${\omega}$-Lie superalgebra. Then ${\rm QC}^{\omega}(g)+[{\rm QC}^{\omega}(g),{\rm QC}^{\omega}(g)]$ is a Lie super ideal of ${\rm GDer}^{\omega}(g)$.
\end{cor}

\begin{proof}Conclusion holds by applying Theorem $3.6$(2),(3).
\end{proof}

\begin{prop}
{\it Let $g$ be
an ${\omega}$-Lie superalgebra. Then $[{\rm C}^{\omega}(g),{\rm
QC}^{\omega}(g)]\subseteq {\rm End}(g,{\rm Z}(g)).$
 Moreover, if ${\rm
Z}(g)=\{0\},$ then $[{\rm C}^{\omega}(g),{\rm QC}^{\omega}(g)]=\{0\}.$}
\end{prop}

\begin{proof}
From similar arguments with theorem in (\cite{Zhang2}, Proposition 2.5), we see that $[{\rm C}(g),{\rm
QC}(g)]\subseteq {\rm End}(g,{\rm Z}(g)).$ So $[{\rm C}^{\omega}(g),{\rm
QC}^{\omega}(g)]\subseteq [{\rm C}(g),{\rm
QC}(g)]\subseteq {\rm End}(g,{\rm Z}(g)).$
\end{proof}

\begin{thm}{\it Let $g$ be
an ${\omega}$-Lie superalgebra. if ${\rm Z}(g)=\{0\}$, then ${\rm C}^{\omega}(g)={\rm QDer}^{\omega}(g)\cap {\rm QC}^{\omega}(g).$}
\end{thm}
\begin{proof}
To prove the first side, we take an arbitrary $d\in{\rm QDer}^{\omega}(g)\cap {\rm QC}^{\omega}(g)$, we recall Proposition $2.6$ in (\cite{Zhang2}, Chapter 2) and assert that $d\in{\rm QDer}(g)\cap {\rm QC}(g)\subseteq{\rm C}(g)$. We see that $d$ is compatible, hence, $d\in{\rm C}^{\omega}(g)$, which means ${\rm QDer}^{\omega}(g)\cap {\rm QC}^{\omega}(g)\subseteq{\rm C}^{\omega}(g).$

By Proposition $2.6$ in (\cite{Zhang2}, Chapter 2), we note that for any $d\in {\rm C}^{\omega}(g)$, $d\in{\rm QDer}(g)\cap{\rm QC}(g).$  It follows that $d\in{\rm QDer}^{\omega}(g)\cap{\rm QC}^{\omega}(g)$ since $d$ is compatible. So ${\rm C}^{\omega}(g)\subseteq{\rm QDer}^{\omega}(g)\cap {\rm QC}^{\omega}(g)$
\end{proof}

\begin{prop}{\it Let $g$ be
an ${\omega}$-Lie superalgebra, if ${\rm Z}(g)=\{0\}$, then ${\rm QC}^{\omega}(g)$ is an Lie superalgebra if and only if $[{\rm QC}^{\omega}(g),{\rm QC}^{\omega}(g)]=0.$}
\end{prop}

\begin{proof} To prove the first side,  we take homogeneous linear maps $d,~d'\in {\rm
QC}^{\omega}(g)$, homogenous element $x\in g.$  Notice that $[d,d']\in {\rm
QC}^{\omega}(g)$, of degree $\mid d\mid\mid d'\mid,$ and
$$[[d,d'](x),y]=[[d,d'](x),y]=(-1)^{(\mid d\mid+\mid d'\mid)\mid x\mid}[x,[d,d'](y)].$$
From the proof of Theorem~3.6~$(3)$, we have
$$[[d,d'](x),y]=[[d,d'](x),y]=-(-1)^{(\mid d\mid+\mid d'\mid)\mid x\mid}[x,[d,d'](y)].$$
Hence
$[[d,d'](x),y]=[[d,d'](x),y]=0,$
i.e. $[d,d']=0.$

The inverse of the inclusion is clear.
\end{proof}

\section{Quasiderivations of $\omega$-Lie
superalgebras}

In this section, we will prove that every compatible quasiderivation of an $\omega$-Lie
superalgebra $g$ can be embedded as a compatible derivations into a larger $\omega$-Lie superalgebra $\tilde{g}$. We also
obtain a direct sum decomposition of Der($g$) when the annihilator
of $g$ is equal to zero.

Let $g=g_{\bar0}\oplus g_{\bar{1}}$ such that $(g,[~,~])$ be an $\omega$-Lie superalgebra over $\mathbb{K}$ and $t$ an indeterminant. We define $\breve{g}:=g[t\mathbb{K}[t]/(t^{3})]=
\{\Sigma(x_{\mu}\otimes t+x_{\lambda}\otimes t^{2}):x_{\mu},x_{\lambda}\in
g\}.$ Then
$\breve{g}$ is an $\omega$-Lie superalgebra with the operation
$[x_{\mu}\otimes t,x_{\lambda}\otimes
t]=[x_{\mu},x_{\lambda}]\otimes t^{2},$ otherwise $[x_{\mu}\otimes t^{i},x_{\lambda}\otimes
t^{j}]=0,$  and $\breve{\omega}(x_{\mu}\otimes t,x_{\lambda}\otimes
t)=\omega(x_{\mu},x_{\lambda}),$ otherwise $\breve{\omega}(x_{\mu}\otimes t^{i},x_{\lambda}\otimes
t^{j})=0$ for all homogeneous elements
$x_{\mu},x_{\lambda}\in g$, $i,j\in\{1,2\}$.

For notational convenience, we write $xt(xt^{2})$ in place of $x\otimes t(x\otimes t^{2}).$

If $U$ is a $\rm Z_{2}$-graded subspace of $g$ such that $g=[g,g]\oplus
U,$ then $\breve g$ has the following decomposition of vector spaces:

$$\breve{g}=gt+gt^{2}=gt+[g,g]t^{2}+Ut^{2},$$ or more
precisely,$$\breve{g}=\breve{g}_{\bar0}\oplus\breve{g}_{\bar1}
=(g_{\bar0}t+[g,g]_{\bar0}t^{2}+U_{\bar0}t^{2})\oplus(g_{\bar1}t+[g,g]_{\bar1}t^{2}+U_{\bar1}t^{2}).$$

Now we define a map $\varphi:{\rm QDer}(g)\rightarrow {\rm
End}(\breve{g})$ satisfying
$$\varphi(d)(at+bt^{2}+ut^{2})=d(a)t+d'(b)t^{2},$$
where homogeneous elements $a\in g, b\in [g,g], u\in U$ have the same degrees, $d$ is a homogenous quasiderivation associated with a linear map $d'$. \vspace{0.3cm}

Immediately, we get the following proposition:

\begin{prop}
$\varphi$ is an even injective map. $\varphi(d)$ is well-defined, and it is independent of the choice of $d'$.
\end{prop}

\begin{prop}
$\varphi(d)\in {\rm Der}(\breve{g}),$ for all $d\in {\rm QDer}(g)$. Moreover, $\varphi(d)$ is compatible $\Leftrightarrow$ $d$ is compatible.
\end{prop}

\begin{proof} From similar argument with theorem in (\cite{Zhang2}, Chapter 4), one gets $\varphi({\rm QDer}(g))\subseteq {\rm Der}(\breve{g}).$

Moreover, we write $(a,b,u)$ for the element $at+bt^{2}+ut^{2}$ for convince. For all $(a,b,u)\in \breve{g}_{\mu}, (a',b',u')\in \breve{g}_{\lambda},$ where $a,a'\in g$, $b,b'\in [g,g],$ $u,u'\in U$, and $a,b,u$ are of the same degree denoted by $\mid a\mid$, $a',b',u'$ are of the same degree denoted by $\mid a'\mid$. $d,d'\in {\rm QDer}(g)$ are homogeneous, it is clear that $\varphi(d)(a,b,u)=(d(a),d'(b),0).$

We observe that
$$\begin{array}{ll}\breve{\omega}(\varphi(d)(a,b,u),(a',b',u'))+(-1)^{\mid d\mid\mid a\mid}\breve{\omega}((a,b,u),\varphi(d)(a',b',u'))\\
=\breve{\omega}((d(a),d'(b),0),(a',b',u'))+(-1)^{\mid d\mid\mid a\mid}\breve{\omega}((a,b,u),(d(a'),d'(b'),0))\\
=\omega(d(a),a')+(-1)^{\mid d\mid\mid a\mid}\omega(a,d(a')).\end{array}$$
 Therefore, $\varphi(d)$ is compatible $\Leftrightarrow$ $d$ is compatible.
\end{proof}

Combining Proposition $4.1$ and Proposition $4.2$, we obtain

\begin{cor}$\varphi({\rm QDer}^{\omega}(g))\subseteq {\rm Der}^{\omega}(\breve{g})$.
\end{cor}

The following theorem describes the precise relationship among  ${\rm Der}^{\omega}(\breve{g})$, ${\rm Zer}(\breve{g})$ and $\varphi({\rm QDer}^{\omega}(g))$.

\begin{thm}Let $(g,[~,~],\omega)$ be an $\omega$-Lie superalgebra with ${\rm Z}(g)=\{0\}$ and $\breve{g},~\varphi$ are as
defined above. Then ${\rm Der}^{\omega}(\breve{g})=\varphi({\rm
QDer}^{\omega}(g))\oplus {\rm ZDer}(\breve{g}).$
\end{thm}
\begin{proof}
The conclusion $\varphi({\rm
QDer}(g))\oplus{\rm ZDer}(\breve{g})\subseteq{\rm Der}(\breve{g})$ follows from similar arguments with (\cite{Zhang2}, Theorem 4.2). By Proposition $3.5$ and Corollary $4.3$, we have $\varphi({\rm
QDer}^{\omega}(g))\oplus{\rm ZDer}(\breve{g})\subseteq{\rm Der}^{\omega}(\breve{g}).$

To prove the converse containment, notice that for any $\tilde d\in {\rm Der}^{\omega}(\breve{g})\subseteq {\rm Der}(\breve{g}),$ by (\cite{Zhang2}, Theorem 4.2),
$\tilde d$ can be written as $\tilde d=\varphi(d')+\tilde f$, where $d'\in {\rm
QDer}(g)$ and $\tilde f\in {\rm ZDer}(\breve{g}).$  Proposition $3.2$ together with Proposition $3.5$ imply that $\varphi(d')$ is compatible, which means $d'\in {\rm
QDer}^{\omega}(g)$ from proposition $4.2$. Hence, ${\rm Der}(\breve{g})\subseteq\varphi({\rm
QDer}^{\omega}(g))\oplus{\rm ZDer}(\breve{g})$.
\end{proof}

\section{Computations and Jordan standard forms in dimension $3$}

It follows from \cite{Zhou6} that any nontrivial complex 3-dimensional
$\omega$-Lie superalgebra must be isomorphic to one the following
algebras: $L_1$, $L_2$, $A_{k}$, $B$, $C_k$, $H$. Computations of $3$-dimensional $\omega$-Lie
algebras($L_1$, $L_2$, $A_{k}$, $B$, $C_k$) have been computed in Chen et al. (\cite{Chen4}, TABLE 1, TABLE 2), so here we need only to consider $\omega$-Lie superalgebra $H$.

Recall that $H$ in (\cite{Zhou6}, Theorem 4.1) satisfies
$[x_1,x_2]=x_1,~[x_1,y]=y,~[x_2,y]=[y,y]=0,~\omega(x_1,x_2)=1,~\omega(y,y)=0$,where $x_1,x_2\in g_{\bar{0}},~y\in g_{\bar{1}},$ $k\in {\mathbb
C}.$

$(1)$ Firstly, we compute the space ${\rm GDer}(H)$. Suppose that
\begin{align*}
d_{\bar0}=\begin{pmatrix}a_{11}&a_{12}&0\cr a_{21}&a_{22}&0\cr 0&0&a_{33}\end{pmatrix};d_{\bar0}'=\begin{pmatrix}a'_{11}&a'_{12}&0\cr a'_{21}&a'_{22}&0\cr 0&0&a'_{33}\end{pmatrix};d_{\bar0}''=\begin{pmatrix}a''_{11}&a''_{12}&0\cr a''_{21}&a''_{22}&0\cr 0&0&a''_{33}\end{pmatrix},
\end{align*}
\begin{align*}
d_{\bar1}=\begin{pmatrix}0&0&b_{13}\cr 0&0&b_{23}\cr b_{31}&b_{32}&0\end{pmatrix};d_{\bar1}'=\begin{pmatrix}0&0&b'_{13}\cr 0&0&b'_{23}\cr b'_{31}&b'_{32}&0\end{pmatrix};d_{\bar1}''=\begin{pmatrix}0&0&b''_{13}\cr 0&0&b''_{23}\cr b''_{31}&b''_{32}&0\end{pmatrix},
\end{align*}

where the action of $d$ on $H$ is given by
\begin{eqnarray*}
  &d_{\bar 0}(x_1)=a_{11}x_1+a_{21}x_2,~&d_{\bar 1}(x_1)=b_{31}y,\\
  &d_{\bar 0}(x_2)=a_{12}x_1+a_{22}x_2 ~&d_{\bar 1}(x_2)=b_{32}y,\\
  &d_{\bar 0}(y)=a_{33}y ~&d_{\bar 1}(y)=b_{13}x_1+b_{23}x_2.
\end{eqnarray*}
The action of $d'_{\bar 0}$, $d'_{\bar 1}$, $d''_{\bar 0}$ and $d''_{\bar 1}$ are defined in the similar way. According to the definition of ${\rm GDer}(H)$,
we obtain
\begin{eqnarray*}
  &a_{11}=a''_{11}-a'_{22}=a''_{33}-a'_{33};a_{22}=a''_{11}-a'_{11};a_{33}=a''_{33}-a'_{11};\\
  &a_{12}=a'_{12}=a''_{21}=0;a_{21}=a'_{21};b_{32}=b'_{32}=b''_{31};b_{31}=b'_{31};\\
  &b_{23}=b'_{23}=b''_{13};b_{13}=b'_{13}=b''_{23}=0.
\end{eqnarray*}

Hence,
\begin{align*}
d=\begin{pmatrix}a_{11}&0&0\cr a_{21}&a_{22}&b_{23}\cr b_{31}&b_{32}&a_{33}\end{pmatrix},
\end{align*}
 together with
\begin{align*}
d'=\begin{pmatrix}a'_{11}&0&0\cr a_{21}&a_{22}-a_{11}+a'_{11}&b_{23}\cr b_{31}&b_{32}&a_{33}-a_{11}+a'_{11}\end{pmatrix}\end{align*}
and \begin{align*}d''=\begin{pmatrix}a_{22}+a'_{11}&a''_{12}&b_{23}\cr 0&a''_{22}&0\cr b_{32}&b''_{32}&a'_{11}+a_{33}\end{pmatrix}
\end{align*}
form a generalized derivation of $H$, which means ${\rm dim}({\rm GDer}(H))=7.$ Since ${\rm dim}({\rm End}(H))=9,$ one gets ${\rm GDer}(H)\neq {\rm End}(H)$.

$(2)$ Secondly, we compute the space ${\rm GDer}^{\omega}(H)$ by applying compatible condition for $\{x_1,x_2,y\}$ and get
$a_{11}+a_{22}=0;b_{23}=b_{13}=0.$ Hence, $d\in {\rm GDer}^{\omega}(H)$ is written as
\begin{align*}
d=\begin{pmatrix}a_{11}&0&0\cr a_{21}&-a_{11}&0\cr b_{31}&b_{32}&a_{33}\end{pmatrix},
\end{align*}
together with
\begin{align*}
d'=\begin{pmatrix}a'_{11}&0&0\cr a_{21}&-2a_{11}+a'_{11}&0\cr b_{31}&b_{32}&a_{33}-a_{11}+a'_{11}\end{pmatrix}\end{align*}
and \begin{align*}d''=\begin{pmatrix}-a_{11}+a'_{11}&a''_{12}&0\cr 0&a''_{22}&0\cr b_{32}&b''_{32}&a'_{11}+a_{33}\end{pmatrix}
\end{align*}

which means ${\rm dim}({\rm GDer}^{\omega}(H))=5.$ So ${\rm GDer}^{\omega}(H)\neq {\rm GDer}(H)$.

$(3)$ We could compute ${\rm QDer}(H)$ by setting $d'=d$ and get that
\begin{eqnarray*}
  &a_{12}=a''_{21}=0;a_{11}+a_{22}=a''_{11};a_{11}+a_{33}=a''_{33};\\
  &b_{13}=b''_{23}=0;b_{32}=b''_{31};b_{23}=b''_{13}.
\end{eqnarray*} So we write a quasiderivation $d$ for
\begin{align*}d=\begin{pmatrix}a_{11}&0&0\cr a_{21}&a_{22}&b_{23}\cr b_{31}&b_{32}&a_{33}\end{pmatrix},~{\rm with}
~d''=\begin{pmatrix}a_{11}+a_{22}&a''_{12}&b_{23}\cr 0&a''_{22}&0\cr b_{32}&b''_{32}&a_{11}+a_{33}\end{pmatrix},\end{align*} which means ${\rm dim}({\rm QDer}(H))=7$, so ${\rm QDer}(H)={\rm GDer}(H)$.

$(4)$ We compute the space ${\rm QDer}^{\omega}(H)$ by adding conditions $a_{11}+a_{22}=0;b_{23}=b_{13}=0$ for $d\in {\rm QDer}(H)$ and have
\begin{align*}d=\begin{pmatrix}a_{11}&0&0\cr a_{21}&-a_{11}&0\cr b_{31}&b_{32}&a_{33}\end{pmatrix}~{\rm with}
~d''=\begin{pmatrix}0&a''_{12}&0\cr 0&a''_{22}&0\cr b_{32}&b''_{32}&a_{11}+a_{33}\end{pmatrix}.\end{align*} So ${\rm dim}({\rm QDer}^{\omega}(H))=5.$ We get that ${\rm QDer}^{\omega}(H)\neq {\rm QDer}(H),$ ${\rm QDer}^{\omega}(H)={\rm GDer}^{\omega}(H).$

From \cite {Chen4}, Corollary 4.6, for all non-Lie 3-dimensional complex $\omega-$Lie algebra $g$, one gets ${\rm GDer}^{\omega}(g)={\rm GDer}(g)$ and ${\rm QDer}^{\omega}(g)={\rm QDer}(g)$, but is not true for the only non-$\omega-$Lie 3-dimensional complex $\omega-$Lie superalgebra $H$.

The rest of this section is devoted to computing the Jordan standard forms of elements in ${\rm GDer}(H)$, ${\rm GDer}^{\omega}(H)$, ${\rm QDer}(H)$ and ${\rm QDer}^{\omega}(H)$.

\begin{prop}
For any $d\in {\rm GDer}(H)$, the Jordan standard form of $d$ is one of the following:
\begin{align*}
(1)\begin{pmatrix}a&0&0\cr 0&b&0\cr 0&0&c\end{pmatrix};~(2)~\begin{pmatrix}a&0&0\cr 1&a&0\cr 0&0&b\end{pmatrix},(3)~\begin{pmatrix}a&0&0\cr 1&a&0\cr 0&1&a\end{pmatrix},~~where~ a,b,c\in \mathbb{C}.
\end{align*}
\end{prop}

\begin{proof}
For any
$d=\begin{pmatrix}a_{11}&0&0\cr a_{21}&a_{22}&b_{23}\cr b_{31}&b_{32}&a_{33}\end{pmatrix}\in {\rm GDer}(H),$  the characteristic polynomial of $d$ is $$f(\lambda)=
(\lambda-a_{11})(\lambda^2-(a_{22}+a_{33})\lambda+a_{22}a_{33}-b_{23}b_{32}).$$ Our arguments will be separated into two cases.

{\bf Case 1} $b_{23}b_{32}=0$. which means $f(\lambda)=(\lambda-a_{11})(\lambda-a_{22})(\lambda-a_{33}).$

{\bf Case 1.1} $a_{11}=a_{22}=a_{33}=a.$ Suppose  $f(\lambda)=(\lambda-a)^{3}$

For $\lambda=a,$ let \begin{align*}B=A-aI=\begin{pmatrix}0&0&0\cr a_{21}&0&b_{23}\cr b_{31}&b_{32}&0\end{pmatrix},\end{align*}
direct calculations show that \begin{align*}B^{2}=\begin{pmatrix}0&0&0\cr b_{23}b_{31}&0&0\cr b_{32}a_{21}&0&0\end{pmatrix},~{\rm and}~B^{3}=0.\end{align*}

{\bf Case 1.1.1}  Suppose $r(B^2)=0$ i.e. $b_{23}b_{31}=b_{32}a_{21}=0,$ combining with $b_{23}b_{32}=0$, one gets
$r(B)=0$ or $1.$ Therefore, \begin{align*}J=\begin{pmatrix}a&0&0\cr 0&a&0\cr 0&0&a\end{pmatrix}~{\rm or}~\begin{pmatrix}a&0&0\cr 0&a&0\cr 0&1&a\end{pmatrix}.\end{align*}

{\bf Case 1.1.2}  Suppose $r(B^2)=1$, which means one of $b_{23}b_{31}$ and $b_{32}a_{21}$ is not $0$. Calculations for any all cases show that $r(B)=2$. So
\begin{align*}J=\begin{pmatrix}a&0&0\cr 1&a&0\cr 0&1&a\end{pmatrix}.\end{align*}

{\bf Case 1.2} Two of $a_{11},a_{22}$ and $a_{33}$ are equal. Suppose $a_{11}=a_{22}=a,$ $a_{33}=b,$ $(a\neq b)$.  Suppose  $f(\lambda)=(\lambda-a)^{2}(\lambda-b)$.

For $\lambda=a,$ let \begin{align*}B=A-aI=\begin{pmatrix}0&0&0\cr a_{21}&0&b_{23}\cr b_{31}&b_{32}&b-a\end{pmatrix},\end{align*}
direct calculations show \begin{align*}B^{2}=\begin{pmatrix}0&0&0\cr b_{23}b_{31}&0&b_{23}(b-a)\cr b_{32}a_{21}+b_{31}(b-a)&b_{32}(b-a)&(b-a)^{2}\end{pmatrix},\end{align*}
\begin{align*}B^{3}=\begin{pmatrix}0&0&0\cr b_{23}b_{31}(b-a)&0&b_{23}(b-a)^{2}\cr b_{32}a_{21}(b-a)+b_{31}(b-a)^{2}&b_{32}(b-a)^{2}&(b-a)^{3}\end{pmatrix},\end{align*}

{\bf Case 1.2.1}  Suppose $r(B)=1$, hence $r(B^2)=r(B^3)=1$, calculations show $r(B^2)+r(I)-2r(B)=2$, so \begin{align*}J=\begin{pmatrix}a&0&0\cr 0&a&0\cr 0&0&b\end{pmatrix}.\end{align*}

{\bf Case 1.2.2}  Suppose $r(B)=2$, hence $b_{23}=0,$ $a_{21}\neq 0,$ which means $r(B^2)=r(B^3)=1$. So \begin{align*}J=\begin{pmatrix}a&0&0\cr 1&a&0\cr 0&0&b\end{pmatrix}.\end{align*}

{\bf Case 1.3} $a_{11},a_{22},a_{33}$ are three different numbers. Suppose $a_{11}=a,$ $a_{22}=b,$ $a_{33}=c,$ $(a\neq b\neq c)$. Therefore, $f(\lambda)=(\lambda-a)(\lambda-b)(\lambda-c)$. In this case, \begin{align*}J=\begin{pmatrix}a&0&0\cr 0&b&0\cr 0&0&c\end{pmatrix}.\end{align*}

{\bf Case 2} $b_{23}b_{32}\neq 0$.

In this case,
$f(\lambda)=(\lambda-a_{11})
(\lambda-\frac{a_{22}+a_{33}\pm\sqrt{(a_{22}-a_{33})^{2}+4b_{23}b_{32}}}{2})$.

{\bf Case 2.1} If $a_{11}$ equals one of $\frac{a_{22}+a_{33}\pm\sqrt{(a_{22}-a_{33})^{2}+4b_{23}b_{32}}}{2}$, take $a_{11}=\frac{a_{22}+a_{33}+\sqrt{(a_{22}-a_{33})^{2}+4b_{23}b_{32}}}{2}=a$, the other one is $b$ for example,
suppose
$f(\lambda)=(\lambda-a)^2(\lambda-b),$ $a\neq b.$

For $\lambda=a,$ let \begin{align*}B=A-aI=\begin{pmatrix}0&0&0\cr a_{21}&a_{22}-a&b_{23}\cr b_{31}&b_{32}&a_{33}-a\end{pmatrix}.\end{align*}

We have \begin{align*}B^{2}=\begin{pmatrix}0&0&0\cr b_{23}b_{31}+(a_{22}-a)a_{21}&(a_{22}-a)^{2}+b_{23}b_{32}&b_{23}(a_{22}+a_{33}-2a)\cr b_{32}a_{21}+b_{31}(a_{33}-a)&b_{32}(a_{22}+a_{33}-2a)&b_{23}b_{32}+(a_{33}-a)^{2}\end{pmatrix}.\end{align*}

{\bf Case 2.1.1}  Suppose $r(B)=1$, hence $r(B^2)=r(B^3)=1$, which means \begin{align*}J=\begin{pmatrix}a&0&0\cr 0&a&0\cr 0&0&b\end{pmatrix}.\end{align*}

{\bf Case 2.1.2}  Suppose $r(B)=2$, we have $r(B^2)=r(B^3)=1$. So \begin{align*}J=\begin{pmatrix}a&0&0\cr 1&a&0\cr 0&0&b\end{pmatrix}.\end{align*}

{\bf Case 2.2} If $a_{11}\neq \frac{a_{22}+a_{33}\pm\sqrt{(a_{22}-a_{33})^{2}+4b_{23}b_{32}}}{2},$ $J=\begin{pmatrix}a&0&0\cr 0&b&0\cr 0&0&c\end{pmatrix}$, where $a\neq b\neq c$.
\end{proof}

\begin{prop}For any $d\in {\rm GDer}^{\omega}(H)$, the Jordan standard form of $d$ is one of the following:
\begin{align*}
(1)\begin{pmatrix}a&0&0\cr 0&-a&0\cr 0&0&b\end{pmatrix};(2)~\begin{pmatrix}a&0&0\cr 1&a&0\cr 0&0&-a\end{pmatrix};(3)~\begin{pmatrix}0&0&0\cr 1&0&0\cr 0&0&c\end{pmatrix};
(4)~\begin{pmatrix}0&0&0\cr 1&0&0\cr 0&1&0\end{pmatrix},
\end{align*}where~ $a,b,c\in \mathbb{C}, c\neq 0.$
\end{prop}

\begin{proof}
For any
$d=\begin{pmatrix}a_{11}&0&0\cr a_{21}&-a_{11}&0\cr b_{31}&b_{32}&a_{33}\end{pmatrix}\in {\rm GDer}^{\omega}(H),$  the characteristic polynomial of $d$ is $f(\lambda)=
(\lambda-a_{11})(\lambda+a_{11})(\lambda-a_{33})$. Our arguments will be separated into two cases.

{\bf Case 1} $a_{11}=-a_{11}=0.$ We get $f(\lambda)=\lambda^{2}(\lambda-a_{33}).$

{\bf Case 1.1} $a_{33}=0$, therefore, $f(\lambda)=\lambda^{3}.$ For $\lambda=0$,
let \begin{align*}B=d=\begin{pmatrix}0&0&0\cr a_{21}&0&0\cr b_{31}&b_{32}&0\end{pmatrix},{\rm then}~~B^{2}=\begin{pmatrix}0&0&0\cr0&0&0\cr a_{21}b_{32}&0&0\end{pmatrix},B^{3}=0.\end{align*}
We obtain that $r(B)=2$$(a_{21}b_{32}\neq 0)$ or $r(B)=1$$(a_{21}b_{32}=0)$ or $r(B)=0$$(d=0)$. So $r(B^2)=1$ or $r(B^2)=0$.

Therefore, $J=\begin{pmatrix}0&0&0\cr 1&0&0\cr 0&1&0\end{pmatrix}$, or $J=\begin{pmatrix}0&0&0\cr 0&0&0\cr 0&1&0\end{pmatrix}$ or $J=0.$

{\bf Case 1.2} $a_{33}\neq0$, so $f(\lambda)=\lambda^{2}(\lambda-a_{33}).$ For $\lambda=0$,
let \begin{align*}B=d=\begin{pmatrix}0&0&0\cr a_{21}&0&0\cr b_{31}&b_{32}&a_{33}\end{pmatrix},\end{align*}
then  \begin{align*}
B^{2}=\begin{pmatrix}0&0&0\cr0&0&0\cr a_{21}b_{32}+a_{33}b_{31}&a_{33}b_{32}&a^{2}_{33}\end{pmatrix}.\end{align*}
We obtain that $r(B)=2$$(a_{21}\neq 0)$ or $r(B)=1$$(a_{21}=0)$, and $r(B^2)=r(B^3)=1$ since $a_{33}\neq0$.

Therefore, $J=\begin{pmatrix}0&0&0\cr 1&0&0\cr 0&0&a_{33}\end{pmatrix}$, or $J=\begin{pmatrix}0&0&0\cr 0&0&0\cr 0&0&a_{33}\end{pmatrix}$.

{\bf Case 2} $a_{11}\neq0.$

{\bf Case 2.1} $a_{33}$ equals to one of  $a_{11}$ and $-a_{11}$. Take $a_{33}=a_{11}$ for instance, $f(\lambda)=(\lambda-a_{11})^{2}(\lambda+a_{11}).$ For $\lambda=a_{11}$,
let \begin{align*}B=d-a_{11}I=\begin{pmatrix}0&0&0\cr a_{21}&-2a_{11}&0\cr b_{31}&b_{32}&0\end{pmatrix},\end{align*}
direct calculations show \begin{align*}B^{2}=\begin{pmatrix}0&0&0\cr -2a_{11}a_{21}&4a^{2}_{11}&0\cr b_{32}a_{21}&-2a_{11}b_{32}&0\end{pmatrix},~{\rm and}~B^{3}=\begin{pmatrix}0&0&0\cr 4a^{2}_{11}a_{21}&-8a^{3}_{11}&0\cr -2b_{32}a_{21}a_{11}&4a^{2}_{11}b_{32}&0\end{pmatrix}.\end{align*}
So $r(B)=2$$(a_{21}b_{32}+2a_{11}b_{31}\neq 0)$ or $r(B)=1$$(a_{21}b_{32}+2a_{11}b_{31}=0)$, and $r(B^2)=1$ since $\mid B^2\mid=0$ and $a_{11}\neq0.$

Therefore, $J=\begin{pmatrix}a_{11}&0&0\cr 1&a_{11}&0\cr 0&0&-a_{11}\end{pmatrix}$, or $J=\begin{pmatrix}a_{11}&0&0\cr 0&a_{11}&0\cr 0&0&-a_{11}\end{pmatrix}$.

The same kind of $J$ could be obtained for case $a_{33}=-a_{11}$, so we omit it.

{\bf Case 2.2} $a_{33}\neq a_{11}\neq-a_{11}$. Therefore, $J=\begin{pmatrix}a_{11}&0&0\cr 0&-a_{11}&0\cr 0&0&a_{33}\end{pmatrix}.$
\end{proof}

For the fact that ${\rm QDer}(H)={\rm GDer}(H)$ and ${\rm QDer}^{\omega}(H)={\rm GDer}^{\omega}(H)$, there are same Jordan standard forms of ${\rm QDer}(H)$(resp. ${\rm QDer}^{\omega}(H)$)and ${\rm QDer}(H)$(resp. ${\rm QDer}^{\omega}(H)$).

For convenience, the characterizations and Jordan standard forms of elements in ${\rm GDer}(H)$, ${\rm GDer}^{\omega}(H)$, ${\rm QDer}(H)$ and ${\rm QDer}^{\omega}(H)$ are summarized in the following Table 1 and Table 2.
\vspace{3mm}

\SMALL{
\hspace*{\fill} Table 1: Characterizations of elements \hspace*{\fill}

\noindent\begin{longtable}{c|c}
 \hline
 $H$& Characterizations of elements\\
  \hline
\tabincell{l}{${\rm dim}({\rm GDer}(H))=7$}&$d=\begin{pmatrix}a_{11}&0&0\cr a_{21}&a_{22}&b_{23}\cr b_{31}&b_{32}&a_{33}\end{pmatrix},$\\
\tabincell{l}{}&${\rm~with}~d'=\begin{pmatrix}a'_{11}&0&0\cr a_{21}&a_{22}-a_{11}+a'_{11}&b_{23}\cr b_{31}&b_{32}&a_{33}-a_{11}+a'_{11}\end{pmatrix}$,\\
\tabincell{l}{}&$~~{\rm and}~~d''=\begin{pmatrix}a_{22}+a'_{11}&a''_{12}&b_{23}\cr 0&a''_{22}&0\cr b_{32}&b''_{32}&a'_{11}+a_{33}\end{pmatrix}$.\\
  \hline
\tabincell{l}{${\rm dim}({\rm GDer}^{\omega}(H))=5$}&$d=\begin{pmatrix}a_{11}&0&0\cr a_{21}&-a_{11}&0\cr b_{31}&b_{32}&a_{33}\end{pmatrix},$\\
\tabincell{l}{}&${\rm~with}~d'=\begin{pmatrix}a'_{11}&0&0\cr a_{21}&-2a_{11}+a'_{11}&0\cr b_{31}&b_{32}&a_{33}-a_{11}+a'_{11}\end{pmatrix}$\\
\tabincell{l}{}&$~~{\rm and}~~d''=\begin{pmatrix}-a_{11}+a'_{11}&a''_{12}&0\cr 0&a''_{22}&0\cr b_{32}&b''_{32}&a'_{11}+a_{33}\end{pmatrix}.$\\
  \hline
\tabincell{l}{${\rm dim}({\rm QDer}(H))=7$}&$d=\begin{pmatrix}a_{11}&0&0\cr a_{21}&a_{22}&b_{23}\cr b_{31}&b_{32}&a_{33}\end{pmatrix},~{\rm with}
~d''=\begin{pmatrix}a_{11}+a_{22}&a''_{12}&b_{23}\cr 0&a''_{22}&0\cr b_{32}&b''_{32}&a_{11}+a_{33}\end{pmatrix},$\\
  \hline
\tabincell{l}{${\rm dim}({\rm QDer}^{\omega}(H))=5$}&$d=\begin{pmatrix}a_{11}&0&0\cr a_{21}&-a_{11}&0\cr b_{31}&b_{32}&a_{33}\end{pmatrix}~{\rm with}
~d''=\begin{pmatrix}0&a''_{12}&0\cr 0&a''_{22}&0\cr b_{32}&b''_{32}&a_{11}+a_{33}\end{pmatrix}.$\\
  \hline
\end{longtable}
}

\SMALL{
\hspace*{\fill} Table 2: Jordan standard forms about elements \hspace*{\fill}

\noindent\begin{longtable}{c|c}
 \hline
 $H$&Jordan standard form about elements \\
  \hline
\tabincell{l}{${\rm GDer}(H)$}&$\begin{pmatrix}a&0&0\cr 0&b&0\cr 0&0&c\end{pmatrix};\begin{pmatrix}a&0&0\cr 1&a&0\cr 0&0&b\end{pmatrix};\begin{pmatrix}a&0&0\cr 1&a&0\cr 0&1&a\end{pmatrix}, a,b,c\in \mathbb{C}$.\\
  \hline
\tabincell{l}{${\rm GDer}^{\omega}(H)$}&$\begin{pmatrix}a&0&0\cr 0&-a&0\cr 0&0&b\end{pmatrix};\begin{pmatrix}a&0&0\cr 1&a&0\cr 0&0&-a\end{pmatrix};\begin{pmatrix}0&0&0\cr 1&0&0\cr 0&0&c\end{pmatrix};\begin{pmatrix}0&0&0\cr 1&0&0\cr 0&1&0\end{pmatrix}, a,b,c\in \mathbb{C}, c\neq 0$.\\
  \hline
\tabincell{l}{${\rm QDer}(H)$}&$\begin{pmatrix}a&0&0\cr 0&b&0\cr 0&0&c\end{pmatrix};\begin{pmatrix}a&0&0\cr 1&a&0\cr 0&0&b\end{pmatrix};\begin{pmatrix}a&0&0\cr 1&a&0\cr 0&1&a\end{pmatrix}, a,b,c\in \mathbb{C}$.\\
  \hline
\tabincell{l}{${\rm QDer}^{\omega}(H)$}&$\begin{pmatrix}a&0&0\cr 0&-a&0\cr 0&0&b\end{pmatrix};\begin{pmatrix}a&0&0\cr 1&a&0\cr 0&0&-a\end{pmatrix};\begin{pmatrix}0&0&0\cr 1&0&0\cr 0&0&c\end{pmatrix};\begin{pmatrix}0&0&0\cr 1&0&0\cr 0&1&0\end{pmatrix}, a,b,c\in \mathbb{C}, c\neq 0$.\\
  \hline
\end{longtable}
}

\newcommand{\noopsort}[1]{} \newcommand{\printfirst}[2]{#1}
  \newcommand{\singleletter}[1]{#1} \newcommand{\switchargs}[2]{#2#1}

\end{document}